\documentclass[14pt,reqno,intlim]{amsart}

\usepackage{amssymb,amsmath}

\usepackage[arrow,matrix,curve]{xy}

\newdir{ >}{{}*!/-5pt/\dir{>}}

\newcommand{\nc}{\newcommand}

\nc{\bC}{\bold{C}} \nc{\bN}{\Bbb{N}} \nc{\cF}{\mathcal{F}}
\nc{\cE}{\mathcal{E}} \nc{\cR}{\mathcal{R}} \nc{\cM}{\mathcal{M}}
\nc{\al}{\alpha} \nc{\bt}{\beta} \nc{\gm}{\gamma} \nc{\dl}{\delta}
\nc{\om}{\omega} \nc{\sg}{\sigma} \nc{\Sg}{\Sigma} \nc{\vf}{\varphi}
\nc{\ve}{\varepsilon} \nc{\os}{\overset} \nc{\ol}{\overline}
\nc{\ul}{\underline} \nc{\us}{\underset} \nc{\sbs}{\subset}
\nc{\bsl}{\backslash} \nc{\Ra}{\Rightarrow}
\nc{\lra}{\longrightarrow} \nc{\all}{\allowdisplaybreaks}

\nc{\Codes}{\operatorname{{\bold{Codes}}}}
\nc{\RegMono}{\operatorname{\mathcal{R}{\rm{eg}\mathcal{M}{\rm{ono}\!}}}}
\nc{\RegEpi}{\operatorname{\mathcal{R}{\rm{eg}\mathcal{E}{\rm{pi}\!}}}}
\nc{\Mn}{\operatorname{\mathcal{M}{\rm{ono}\!}}}
\nc{\Ep}{\operatorname{\mathcal{E}{\rm{pi}\!}}}
\nc{\Rg}{\operatorname{\mathcal{R}{\rm{eg}\!}}}
\nc{\Ob}{\operatorname{Ob\!}}

\numberwithin{equation}{section}

\newtheorem{theo}{\ \ \ Theorem}[section]

\newtheorem{prop}[theo]{\ \ \ Proposition}

\theoremstyle{definition}

\theoremstyle{remark}
\newtheorem{rem}[theo]{\ \ \ Remark}

\begin{document}

\title[]
{Idempotent radicals and factorization systems}

\author{Dali Zangurashvili}

\maketitle

\begin{abstract}
Idempotent (hereditary/cohereditary)
(pre)radicals on the category of left modules over an associative ring with unity are described in terms of factorization systems. For radicals of the aforementioned kind, this is done in two different ways.

\vskip+2mm
\noindent{\bf Key words and phrases}:  idempotent radical; factorization system; closure operator; left module over an associative ring with identity.
\vskip+2mm

\end{abstract}

\vskip+5mm

\section{Introduction}
The aim of the present paper is to give the description
of idempotent (hereditary/cohereditary) (pre)radicals on the category of left modules over an associative ring with unity in terms of
factorization systems. To this end, the following approach is employed.  

The two injective mappings (call them $F_1$ and $F_2$) from the class of idempotent radicals to the class of idempotent weakly hereditary closure operators are well known (see, e.g., \cite{K} and \cite{DT}). They send an idempotent radical $\mathbf{R}$ to the closure operators $c_{\mathbf{R}}$ and $c^{\mathbf{R}}$, respectively, defined as follows:
$$(c_{\mathbf{R}})_M(N)=N+\mathbf{R}(M),$$
and 
$$(c^{\mathbf{R}})_M(N)=\pi^{-1}(\mathbf{R}(M/N)),$$
for any submodule $N$ of a module $M$ ($\pi:M\twoheadrightarrow M/N$ is the canonical projection).

Moreover, there is a bijection $G$ from the class of idempotent weakly hereditary closure
operators on the category $\Lambda$-$Mod$ to the class of its  factorization
systems with the property that any homomorphism in the right-hand class is a  monomorphism \cite{DT}. Therefore, the compositions of  $F_1$ and $F_2$ with $G$ induce
the bijections between the class of idempotent radicals and  certain
classes of factorization systems. In the present paper, conditions
on factorization systems under which they belong to the latter classes are found. Conditions under
which factorization systems correspond to  hereditary/cohereditary
idempotent radicals are also found.

Note that the factorization system corresponding to an idempotent radical
under the composition $GF_1$ coincides with the one constructed
in \cite{Z3}, \cite{Z4} with the aid of completely different techniques (namely,
those of transporting factorization systems from one category to another via an adjunction,
developed in \cite{Z1}-\cite{Z3}).

Finally note that $c_{\mathbf{R}}$ is an idempotent weakly hereditary closure operator even in the case where $\mathbf{R}$ is an arbitrary idempotent \textit{pre}radical \cite{DT}. Employing this fact, the description of  idempotent/hereditary preradicals in terms of factorization systems are also provided.

Financial support from Shota Rustaveli Georgian National Science
Foundation (Ref.: FR-24-8249) is gratefully acknowledged.


\section{Preliminaries}

We use the terms '(pre)radical', `factorization system', and `closure operator' in the sense of \cite{M},  \cite{AHS}, and \cite{DT} respectively. They are defined in the context of an arbitrary category (satisfying certain natural requirements). But, in this paper, we deal with only the category $\Lambda$-$Mod$ of left $\Lambda$-modules over an associative ring $\Lambda$ with unity. So, we give the definitions only for this particular case.

The facts provided in this section can be found in the aforementioned works and also in \cite{K}.

Throughout the paper,  a module means a left unital module over $\Lambda$. The symbol $Mono$ denotes the class of all monomorphisms of the category $\Lambda$-$Mod$.
\vskip+1mm
1. \textit{Radicals}. A preradical on the category $\Lambda$-$Mod$ is defined as an endofunctor $\mathbf{R}$ of $\Lambda$-$Mod$
such that $\mathbf{R}(M)$ is a submodule of $M$,  for any module $M$, and the embedding $\mathbf{R}(M)\rightarrowtail M$ is natural. 
\vskip+1mm
If $\mathbf{R}(\mathbf{R}(M))=M$, for any module $M$, then $\mathbf{R}$ is called idempotent. 
 
 \vskip+1mm
A preradical $\mathbf{R}$ is called hereditary (resp. cohereditary) if 
\begin{equation}
\mathbf{R}(N)=\mathbf{R}(M)\cap N,
\end{equation}
(resp.
\begin{equation}
\mathbf{R}(M/N)=(N+\mathbf{R}(M))/N,)
\end{equation}
for any submodule $N$ of a module $M$.
\vskip+1mm

Any hereditary preradical is idempotent.

A preradical is called a radical if $$\mathbf{R}(M/\mathbf{R}(M))=0,$$ for any module $M$.

Any cohereditary preradical is a radical.

\vskip+1mm
2. \textit{Closure operators}. Recall that a closure operator on the category $\Lambda$-$Mod$ is a family  $c=(c_M)_{M\in \Lambda-Mod}$ of mappings $c_M$ from the set of submodules of a module $M$ to itself such that the following conditions are satisfied:

-- $N\subseteq c_M(N)$, for any submodule $N$ of a module $M$;

-- if $N\subseteq N'$, then $c_M(N)\subseteq c_M(N')$, for any submodules $N$ and $N'$ of a module $M$;

-- for any homomorphism $f:M\rightarrow M'$ and any submodule $N$ of $M$, $f(c_M(N))\subseteq c_{M'}(f(N))$.
\vskip+1mm
A monomorphism $f:M\rightarrowtail M'$ is called $c$-closed if $c_{M'}(f(M))=f(M)$.

 A homomorphism $f:M\rightarrow M'$ is called $c$-dense if $c_{M'}(f(M))=M'$. 

A closure operation $c$ is called idempotent if $c_M(c_M(N))=c_M(N)$, for any submodule $N$ of a module $M$.

A closure operation $c$ is called weakly hereditary if the inclusion $j_N:N\rightarrowtail c_M(N)$ is dense, for any submodule $N$ of a module $M$.  

A closure operation $c$ is called additive if one has $c_M(N_1+N_2)=c_M(N_1)+c_M(N_2)$, for any submodules $N_1$ and $N_2$ of a module $M$.
\vskip+1mm
3. \textit{Factorization systems}. A factorization system on the category $\Lambda$-$Mod$ is a pair $(\mathbb{E},\mathbb{M})$ of classes of homomorphisms such that:

-- the classes $\mathbb{E}$ and $\mathbb{M}$ are closed under composition with isomorphisms; 

--  every homomorphism
$f$ admits an $(\mathbb{E},\mathbb{M})$-factorization, i.e., there are homomorphisms $e\in \mathbb{E}$ and $m\in \mathbb{M}$ with $f=me$;

-- the diagonalization condition holds,
i.e., for any $e\in \mathbb{E}$ and $m\in \mathbb{M}$, one has $e\downarrow m$ (or, equivalently, $m\uparrow e$),
which means that for any commutative diagram
\begin{equation}
\xymatrix{K\ar[r]^{e}\ar[d]&L\ar[d]\ar@{-->}[dl]^{\delta}\\
N\ar[r]_{m}&M}
\end{equation}
there exists a unique homomorphism $\delta$ making the diagram commutative.
\vskip+1mm

We define an $\mathbb{M}$-submodule of a module $M$ as a submodule $N$ such that its inclusion monomorphism  $N\rightarrowtail M$ lies in $\mathbb{M}$.
\vskip+1mm

Recall some properties of factorization systems.
\vskip+1mm

(a) for any homomorphism $f$, an $(\mathbb{E},\mathbb{M})$-factorization is unique up to isomorphism, i.e., if $f=me=m'e'$, then there exists a unique isomorphism $\delta$ such that the following diagram is commutative:
\begin{equation}
\xymatrix{K\ar[r]^{e}\ar[d]_{e'}&L\ar[d]^{m}\ar@{-->}[dl]^{\delta}_{\simeq}\\
N\ar[r]_{m'}&M}
\end{equation}

(b) one has the equalities $\mathbb{E}=\mathbb{M}^{\uparrow}$ and $\mathbb{M}=\mathbb{E}^{\downarrow}$, where $\mathbb{M}^{\uparrow}$ denotes the class of all homomorphisms $f$ with $f\downarrow m$, for all $m\in \mathbb{M}$; the class 
$\mathbb{M}^{\downarrow}$ is defined dually.
\vskip+1mm
(c) the intersection $\mathbb{E}\cap \mathbb{M}$ is equal to the class of isomorphisms;

(d) the classes $\mathbb{E}$ and $\mathbb{M}$ are closed under composition;

(e)  the class $\mathbb{M}$ is stable under pullback, i.e., if a square
\begin{equation}
\xymatrix{K\ar[r]^{m'}\ar[d]_{f}&L\ar[d]^{g}\\
N\ar[r]^{m}&M}
\end{equation}
is a pullback and $m\in \mathbb{M}$, then $m'\in \mathbb{M}$. In particular, the inverse image of any $\mathbb{M}$-submodule under any homomorphism is an $\mathbb{M}$-submodule. If $\mathbb{M}\subseteq Mono$ , then the latter condition implies the property (e).

The class $\mathbb{E}$ satisfies the property dual to (e);

(f) the class $\mathbb{M}$ satisfies the following cancellation property: for any  homomorphisms $f:C\rightarrow C'$ and $g:C'\rightarrow C''$, if a homomorphism $g$ and the composition $gf$ lie in $\mathbb{M}$, then so is $f$. 
The class $\mathbb{E}$ satisfies the dual property;

(g) the intersection of any set of $\mathbb{M}$-submodules is an $\mathbb{M}$-submodule. In particular, any module has the smallest $\mathbb{M}$-submodule. 
The class $\mathbb{E}$ satisfies the dual property.


\vskip+1mm
A partial order on the class of factorization systems is introduced as follows: $(\mathbb{E},\mathbb{M})\leq (\mathbb{E'},\mathbb{M'})$ if $\mathbb{M}\leq \mathbb{M'}$.
\vskip+2mm

4. \textit{Closure operators and factorization systems.} Any idempotent weakly hereditary closure operator $c$ on $\Lambda$-$Mod$ induces a certain factorization system $G(c)=(\mathbb{E},\mathbb{M})$ on this category. The class $\mathbb{E}$ consists of $c$-dense homomorphisms; while the class $\mathbb{M}$ consists of $c$-closed monomorphisms. 

The $(\mathbb{E},\mathbb{M})$-factorization of a homomorphism $f:M\rightarrow M'$ is given by $f=(ej_{Im\;f})i$, where $f=em$ is the epi-mono factorization of $f$, while $j_{im\;f}$ is the inclusion $Im\;f\rightarrowtail c_B(Im\;f)$:
\begin{equation}
\xymatrix{M\ar[rr]^{f}\ar@{->>}[rd]^{e}&&M'\\
&Im \;f\ar@{ >->}[r]^{j_{Im\;f}}\ar@{ >->}[ur]^{m}&c_B(Im\;f)\ar@{ >->}[u]^{i}}
\end{equation}

 Conversely, any factorization system $(\mathbb{E},\mathbb{M})$ on  $\Lambda$-$Mod$ with $\mathbb{M}\subseteq Mono$ induces an idempotent weakly hereditary closure operator $c$. It is defined as follows: for a submodule $N$ of a module $M$, $c_M(N)$ is the domain of the monomorphism $m$ in the $(\mathbb{E},\mathbb{M})$-factorization of the embedding $N\rightarrowtail M$.
 
 The constructed mappings determine mutually inverse bijections between the class of idempotent weakly hereditary closure operators and that of factorization systems $(\mathbb{E},\mathbb{M})$ with $\mathbb{M}\subseteq Mono$.

\section{Idempotent radicals and factorization systems, I}


As  noted in the Introduction, $c_{\mathbf{R}}$ is an idempotent weakly hereditary closure operator if $\mathbf{R}$ is (not necessarily an idempotent radical, but) an arbitrary idempotent preradical \cite[Proposition 3.4]{DT}. The homomorphisms from the corresponding  factorization system $(\mathbb{E}_1,\mathbb{M}_1)=G(c_{\mathbf{R}})$ can be characterized as follows:
 \vskip+1mm
 a homomorphism $e:M\rightarrow M'$ lies in $\mathbb{E}_1$ if and only if $Im\; e+\mathbf{R}(M')=M'$; 
 
 a homomorphism $m:M\rightarrow M'$ lies in $\mathbb{M}_1$ if and only if it is a monomorphism and $\mathbf{R}(M')\subseteq Im M$. 
 \vskip+1mm
 
 The $(\mathbb{E}_1,\mathbb{M}_1)$-factorization of a homomorphism $f:M\rightarrow M'$ is given by $f=me$ with the homomorphisms $e$ and $m$ indicated in the following diagram:
\begin{equation}
\xymatrix{M\ar[rr]^{f}\ar[rd]^{e}&&M'\\
&Im\;f+\mathbf{R}(M')\ar@{ >->}[ur]^{m}}
\end{equation}
The fact that $e\in \mathbb{E}$ follows from the following inclusions: $$\mathbf{R}(Im\; f+\mathbf{R}(M'))\subseteq \mathbf{R}(M')$$ and  $$\mathbf{R}(M')=\mathbf{R}\mathbf{R}(M')\subseteq \mathbf{R}(Im \;f+\mathbf{R}(M')),$$ whence \begin{equation}\mathbf{R}(Im\;f+\mathbf{R}(M'))=\mathbf{R}(M').
\end{equation} 
\vskip+1mm

Note that $\mathbf{R}(M)$ is the smallest $\mathbb{M}$-submodule of a module $M$. 
 
\vskip+1mm
Consider now an arbitrary factorization system $(\mathbb{E},\mathbb{M})$, and the following conditions on it:
\vskip+1mm
(i)  $\mathbb{M}\subseteq Mono$;
\vskip+1mm
(ii) the class $\mathbb{M}$ satisfies the following cancellation property: if $m=m'f$ with $m\in \mathbb{M}$ and $m'$ being a monomorphism, then $m'\in \mathbb{M}$;
\vskip+1mm
(iii) for any module $M$, the zero homomorphism $0\rightarrowtail M/N$ lies in $\mathbb{M}$, where $N$ is the  smallest $\mathbb{M}$-submodule of $M$.
\vskip+3mm

Note that the conditions (i) and (ii) are satisfied if and only if the condition (i) and the following one are satisfied:
\vskip+1mm
(ii') for any module $M$, a submodule of $M$ is an $\mathbb{M}$-submodule if and only if it contains the smallest $\mathbb{M}$-submodule.

\begin{theo} \label{theo3.1}
(a) There is an order-reversing bijection between the class of idempotent preradicals and the class of factorization systems satisfying the conditions (i)-(ii).  

(b) There is an order-reversing bijection between the class of idempotent radicals and the class of factorization systems satisfying the conditions (i)-(iii).

\end{theo}

\begin{proof}

 Let  $(\mathbb{E},\mathbb{M})$ be a factorization system satisfying the conditions (i)-(ii). For a module $M$, define $\mathbf{R}(M)$ as the smallest $\mathbb{M}$-submodule of $M$. It exists due to the property (g) of factorization systems. To show the functoriality of $\mathbf{R}$, consider a homomorphism $f:M\rightarrow M'$. According to the property (e), the inverse image of  any  $\mathbb{M}$-submodule $N'$ of $M'$ lies in  $\mathbb{M}$. This implies that $f(\mathbf{R}(M))\subseteq \mathbf{R}(M')$.

The idempotency of $\mathbf{R}$ follows from the property (d) of factorization systems. 

According to the property (f), the embedding $\mathbf{R}(M)\rightarrowtail M$ lies in $\mathbb{M}$. Therefore, the condition (ii) implies that the class  $\mathbb{M}_1$ of homomorphisms in the factorization system $(\mathbb{E}_1,\mathbb{M}_1)=GF_1(\mathbf{R})$  is $\mathbb{M}$. Since, in a factorization system, the right-hand class uniquely determines the left-hand class (as it follows from the property (b)), we obtain that $(\mathbb{E}_1,\mathbb{M}_1)=(\mathbb{E},\mathbb{M})$.
\vskip+1mm

If (iii) is satisfied, then  obviously, for any module $M$,  $\mathbf{R}(M/\mathbf{R}(M))=0$.

The fact that if $\mathbf{R}$ is an idempotent preradical (resp. radical), then $(\mathbb{E}_1,\mathbb{M}_1)$ satisfies the conditions (i) and (ii) (resp. (i)-(iii)) is obvious. 
\vskip+1mm

\end{proof}

\begin{rem} Theorem 3.1(a) and Theorem 3.1(b) imply that the condition (iii) does not follow from the conditions (i) and (ii).
\end{rem}

Consider again an arbitrary factorization system $(\mathbb{E},\mathbb{M})$, and the following conditions on it:
\vskip+1mm
(iv) for any monomorphism $f:M\rightarrowtail M'$, and the smallest $\mathbb{M}$-submodule $N'$ of $M'$, the inverse image of $N'$ under $f$ is the smallest $\mathbb{M}$-submodule of $M$;
\vskip+1mm
(v) for any monomorphism $f:M\rightarrowtail M'$, the lower inclined quadrangle in the following commutative diagram is a pullback:
\begin{equation}
\xymatrix{0\ar[rr]^{=}\ar@{ >->}[dd]^{i}\ar@{ >->}[dr]^{e}&&0\ar@{ >->}[dd]^{i'}\ar@{ >->}[dr]^{e'}&\\
&X\ar@{-->}[rr]^{\delta\;\;\;\;\;\;\;\;\;\;\;\;}\ar@{ >->}[dl]^{m}&&Y\ar@{ >->}[dl]^{m'}\\
M\ar@{ >->}[rr]^{f}&&M'}
\end{equation}
Here $me$ and $m'e'$ are $(\mathbb{E},\mathbb{M})$-factorization of the embeddings $i$ and $i'$ of $0$ into $N$ and $M$ respectively.

\begin{prop} \label{prop3.2}
For factorization systems  $(\mathbb{E},\mathbb{M})$ with  $\mathbb{M}\subseteq Mono$, the conditions (iv) and (v) are equivalent. An idempotent preradical $\mathbf{R}$  is hereditary if and only if its factorization system $(\mathbb{E}_1,\mathbb{M}_1)$ satisfies these conditions.
\vskip+1mm
 
\end{prop}

\begin{proof}
The equivalence of (iv) and (v) follows from the fact that, if $\mathbb{M}\subseteq Mono$, then the smallest $\mathbb{M}$-submodule of a module $M$ is the submodule $X$ in the $(\mathbb{E},\mathbb{M})$-factorization of the embedding $i:0\rightarrowtail M$:
\begin{equation}
\xymatrix{0\ar@{ >->}[rr]^{i}\ar@{ >->}[rd]_{e}&&M\\
&X\ar@{ >->}[ur]_{m}&}
\end{equation}
Indeed, if $N$ is an $\mathbb{M}$-submodule of $M$, then there is a monomorphism $\delta$ in the commutative diagram
\begin{equation}
\xymatrix{0\ar[ddr]\ar@{ >->}[rr]^{i}\ar@{ >->}[rd]^{e}&&M\\
&X\ar@{-->}[d]^{\delta}\ar@{ >->}[ur]^{m}&\\
&N\ar@{ >->}[uur]&}
\end{equation}
due to the diagonalization property of a factorization system.
\end{proof}

Theorem \ref{theo3.1} and Proposition \ref{prop3.2} imply the following statement.

\begin{theo}
 The class of  hereditary preradicals (resp. radicals) is bijective to the class of factorization systems satisfying the conditions (i)-(ii)(resp. (i)-(iii)) and the equivalent conditions (iv), (v).

\end{theo}

Further, we introduce yet another condition on a factorization system $(\mathbb{E},\mathbb{M})$:
\vskip+1mm
(vi) for any monomorphism $f:M\rightarrowtail M'$ and its $(\mathbb{E},\mathbb{M}$)-factorization $f=me$, the unique homomorphism $h$ making the following diagram commutative lies in $\mathbb{M}$:
\begin{equation}
\xymatrix{&A\ar@{->>}[r]^{\pi'}\ar@{ >->}[d]^{m}&A/M\ar@{-->}[d]^{h}\\
M\ar@{ >->}[r]_{f}\ar@{ >->}[ur]^{e}&M'\ar@{->>}[r]_{\pi}&M'/M}
\end{equation}

\begin{prop}\label{prop3.3}
An idempotent radical $\mathbf{R}$ is cohereditary if and only if the factorization system $(\mathbb{E}_1,\mathbb{M}_1)$ satisfies the condition (vi).
\vskip+1mm

\end{prop}

\begin{proof} 
According to the construction of the $(\mathbb{E}_1,\mathbb{M}_1)$-factorizations given at the beginning of this section, one has:
\begin{equation}
A\simeq M+\mathbf{R}(M').
\end{equation}
\vskip+1mm 
\noindent This implies "only if" part of the statement. 

For "if" part, consider any monomorphism $f: M\rightarrowtail M'$ and the diagram 
\begin{equation}
\xymatrix{\mathbf{R}(M')\ar@{ >->}[dd]_{i_{\mathbf{R}(M')}}\ar@{ >-->}[dr]^{(\mu_{M'},\mathbf{R}(\pi))}\ar@{ >->}[dddr]\ar[drrr]^{\mathbf{R}(\pi)}&&&\\
&P\ar@{->>}[rr]^{\beta}\ar@{ >->}[dd]&&\mathbf{R}(M'/M)\ar@{ >->>}[dd]^{\mu_{M'/M}}\\
M+\mathbf{R}(M')\ar@{->>}[rr]^{\;\;\;\;\;\;\;\;\;\;\;\pi'}\ar@{ >->}[dr]_{\;\;\;\;\;\;\;\;\;m}\ar@{ >-->}[ur]^{\;\;\;\;\;\;\;\;\;\;\;\;\;\;\;\;\;\;\;\alpha}&&(M+\mathbf{R}(M'))/M\ar@{ >-->}[ru]^{\gamma}\ar@{ >->}[rd]^{h}&\\
M\ar@{ >->}[u]^{i_M}\ar@{ >-->}[uur]_{(f,0)}\ar@{ >->}[r]_{f}&M'\ar@{->>}[rr]_{\pi}&&M'/M}
\end{equation}
where the internal rectangle is a pullback, $\mu_{M'}$ is the embedding $\mathbf{R}(M')\rightarrowtail M'$. The monomorphism  $\gamma$ exists due to the fact that $\mu_{\mathbf{R}(M'/M)}\beta\alpha i_M=0$ and hence $\beta\alpha i_M=0$. Therefore, $(M+\mathbf{R}(M'))/M\subseteq \mathbf{R}(M'/M)$.

For the converse inclusion, recall that the embedding $\mathbf{R}(M'/M)$ is the smallest $\mathbb{M}$-subobject of $M'/M$. Moreover, the existence of isomorphism (3.4) and the condition (v) imply that the homomorphism $h$ lies in $\mathbb{M}$. This implies the claim.
\end{proof}

Theorem 3.1 and Proposition 3.5 imply the following statement.

\begin{theo}
The class of idempotent cohereditary radicals  is bijective to the class of factorization systems satisfying the conditions (i)-(ii) and (vi).
\end{theo}

\begin{proof}
It suffices to recall that any cohereditary preradical is a radical.
\end{proof}

\section{Idempotent radicals and factorization systems, II}

In this section, we deal with the mapping $GF_2$.
 It sends an idempotent radical $\mathbf{R}$ to the factorization system $(\mathbb{E}_2,\mathbb{M}_2)=G(c^{\mathbf{R}})$, where:
 \vskip+1mm
 a homomorphism $e:M\rightarrow M'$ lies in $\mathbb{E}_2$ if and only if $\mathbf{R}(M'/Im\; e)=M'/Im\; e$; 
 
a homomorphism $m:M\rightarrow M'$ lies in $\mathbb{M}_2$ if and only if it is a monomorphism and $\mathbf{R}(M'/Im\; m)=0$. 
 \vskip+1mm
 
 The $(\mathbb{E}_2,\mathbb{M}_2)$-factorization of a homomorphism $f:M\rightarrow M'$ is given by $f=me$:
\begin{equation}
\xymatrix{&P\ar@{->>}[r]^{\pi'\;\;\;\;\;\;\;\;}\ar@{ >->}[d]^{m}&\mathbf{R}(M'/Im\; f)\ar@{ >->}[d]\\
M\ar[r]_{f}\ar@{ -->}[ur]^{e}&M'\ar@{->>}[r]_{\pi\;\;\;\;\;\;\;\;}&M'/Im\; f}
\end{equation} 
In this diagram, the internal square is a pullback, and $e$ is the unique homomorphism with $\pi'e=0$. 

 \vskip+1mm
 
One can easily see that  $(\mathbb{E}_2,\mathbb{M}_2)\leq (\mathbb{E}_1,\mathbb{M}_1)$.
\vskip+1mm
Before continue, recall the following theorem.

\begin{theo} \cite[Theorem 7.2]{K}
An idempotent closure operator $c$ has the form $c^{\mathbf{R}}$ with an idempotent radical $\mathbf{R}$, if and only if it is weakly hereditary and satisfies the following condition: if $K$ and $N$ are submodules of a module $M$ with $K\subseteq N$, then $$c_{M/K}(N/K)=c_{M}(N)/K.$$

\end{theo}
\vskip+2mm

In view of Theorem 4.1, we introduce the following condition on a factorization system $(\mathbb{E},\mathbb{M})$:
\vskip+1mm
(vii) for any submodules $K$ and $N$ of a module $M$ with $K\subseteq N$, the unique homomorphism $\delta$ making the following diagram commutative is the cokernel of the composition $ei$:


\begin{equation}
\xymatrix{&N\ar@{->>}[rr]\ar@{ >->}[dd]\ar@{ >->}[rd]^{e}&&N/K\ar@{ >->}[rd]^{e'}\ar@{ >->}[dd]&\\
&&X\ar@{-->}[rr]^{\;\;\;\;\;\delta}\ar@{ >->}[ld]^{m}&&Y\ar@{ >->}[ld]^{m'}\\
K\ar@{ >->}[r]\ar@{ >->}[uur]^{i}&M\ar@{->>}[rr]&&M/K}
\end{equation}
Here $i$ is the inclusion monomorphism $K\rightarrowtail N$, while $me$ and $m'e'$ are the $(\mathbb{E},\mathbb{M})$-factorizations of the inclusions $N\rightarrowtail M$ and $N/K\rightarrowtail M/K$ respectively. The homomorphism $\delta$ exists due to the diagonalization property of $(\mathbb{E},\mathbb{M})$.

\vskip+2mm
Theorem 4.1 implies the following statement.

\begin{theo}
There exists an order-reversing bijection between the class of idempotent radicals and the class of factorization systems satisfying the conditions (i) and (vii):

\vskip+1mm

\end{theo}

We will now deal with the issue which factorization systems correspond to hereditary radicals. To this end, recall the following proposition.

\begin{prop} \cite[Proposition 7.8]{K}
Let $\mathbf{R}$ be an idempotent radical. Then the following conditions are equivalent:
\vskip+1mm
(a) $\mathbf{R}$ is hereditary;
\vskip+1mm
(b) if $K$ and $N$ are submodules of a module $M$ with $K\subseteq N$ and  $N\subseteq c^{\mathbf{R}}_M(K)$, then $c^{\mathbf{R}}_N(K)=N$;
\vskip+1mm
(c) if $K$ and $N$ are submodules of a module $M$ with $K\subseteq N$, then $c^{\mathbf{R}}_N(K)=N\cap c^{\mathbf{R}}_M(K)$;
\vskip+1mm
(d) for any module $M$, any submodule $N$ of $M$, and any $a\in M$, one has the equality 
$$c^{\mathbf{R}}_{\Lambda}(N:a)=(c^{\mathbf{R}}_{M}(N):a),$$
where  $(N:a)$ denotes the left ideal $\lbrace \lambda\in \Lambda\vert \lambda a\in N\rbrace$ of $\Lambda$.
\end{prop}

Consider an arbitrary factorization system $(\mathbb{E},\mathbb{M})$, and the following conditions on it:
\vskip+1mm
(viii) the class $\mathbb{E}\cap Mono$ satisfies the following cancellation property: if $e=fe'$, $e\in \mathbb{E}$, and $f$ and $e'$ are  monomorphisms, then $e'\in \mathbb{E}$.
\vskip+1mm
(ix) for any submodules $K$ and $N$ of a module $M$ with $K\subseteq N$, and the $(\mathbb{E},\mathbb{M})$-factorization $me$ of the inclusion $K\rightarrowtail M$, the  unique homomorphism $\varphi$ making the following diagram commutative lies in $\mathbb{E}$:
\begin{equation}
\xymatrix{K\ar@{ >->}[rrd]^{e}\ar@{-->}[dr]^{\varphi}\ar@{ >->}[rdd]&&\\
&P\ar@{ >->}[r]\ar@{ >->}[d]&X\ar@{ >->}[d]^{m}\\
&N\ar@{ >->}[r]&M}
\end{equation} Here the internal square is a pullback.

\vskip+1mm

(x) for any module $M$ and any element $a$ in it, the order-preserving mapping $$N\mapsto (N:a)$$ from the set of submodules of $M$ to the set of left ideals of the ring $\Lambda$, preserves $(\mathbb{E},\mathbb{M})$-factorizations in the following sense: if 
\begin{equation}
\xymatrix{N\ar@{ >->}[rr]\ar@{ >->}[dr]_{e}&&M\\
&X\ar@{ >->}[ur]_{m}&}
\end{equation}
is the  $(\mathbb{E},\mathbb{M})$-factorization of the embedding of a submodule $N$ into $M$, then
 \begin{equation}
\xymatrix{(N:a)\ar@{  >->}[rr]\ar@{ >->}[dr]_{e'}&&\Lambda\\
&(X:a)\ar@{ >->}[ur]_{m'}&}
\end{equation}
is the  $(\mathbb{E},\mathbb{M})$-factorization of the embedding $(N:a)\rightarrowtail \Lambda$.

\begin{prop}
The conditions (viii)-(x) are equivalent to one another. Moreover, an idempotent radical $\mathbf{R}$  is hereditary if and only if the factorization system $(\mathbb{E}_2,\mathbb{M}_2)$ satisfies these equivalent conditions.

\end{prop}

\begin{proof}
The equivalence of the conditions (viii)-(x) follows from Proposition 4.3.

The fact that the condition (viii) implies that $\mathbf{R}$  is hereditary follows from the implication (b)$\Rightarrow$(a) in Proposition 4.3. For the converse, assume that $\mathbf{R}$  is hereditary, and consider monomorphisms $e':K\rightarrowtail N$ and $f:N\rightarrowtail M$ with $e=fe'\in \mathbb{E}$. For the closure operator $c^{\mathbf{R}}=G^{-1}(\mathbb{E}_2,\mathbb{M}_2)$, one has $c^{\mathbf{R}}_M(K)=M$. The condition (b) of Proposition 4.3 implies that  $c^{\mathbf{R}}_N(K)=N$, and hence $e'\in \mathbb{E}$. 
\end{proof}

Theorem 4.2 and Proposition 4.4 imply the following statement.

\begin{theo}
The class of  hereditary radicals is bijective to the class of factorization systems satisfying the conditions (i), (vii) and the equivalent conditions (viii)-(x).
\end{theo}

Before continue, recall the following theorem.

\begin{theo} \cite[Theorem 3.4 (6)]{DT} Let $\mathbf{R}$ be an idempotent radical. Then the following conditions are equivalent:
\vskip+1mm
(a) $\mathbf{R}$ is cohereditary;
\vskip+1mm
(b) $c^{\mathbf{R}}$ is additive;
\vskip+1mm
(c) $c_{\mathbf{R}}=c^{\mathbf{R}}$.
\end{theo}

In view of the condition (b) of Theorem 4.6, we introduce the following condition on a factorization system $(\mathbb{E},\mathbb{M})$:
\vskip+1mm
(xi) for any module $M$, any submodules $N_1$ and $N_2$ of $M$, and the $(\mathbb{E},\mathbb{M})$-factorizations $m_1e_1$ and $m_2e_2$ of their inclusion monomorphisms, one has $i\in \mathbb{E}$ and $i'\in \mathbb{M}$, where $i$ and $i'$ are the homomorphisms indicated in the commutative diagram
\begin{equation}
\xymatrix{N_1+N_2\ar@{ >->}[dr]^{i}&&N_2\ar@{ >->}[d]^{e_2}\ar@{ >->}[ll]\\
&X_1+X_2\ar@{ >->}[dr]^{i'}&X_2\ar@{ >->}[l]\ar@{ >->}[d]^{m_2}\\
N_1\ar@{ >->}[uu]\ar@{ >->}[r]_{e_1}&X_1\ar@{ >->}[r]_{m_1}\ar@{ >->}[u]&M}
\end{equation}
\vskip+2mm

\begin{prop}
An idempotent radical $\mathbf{R}$ is cohereditary if and only if the following equivalent conditions are satisfied:
\vskip+1mm
(a) the factorization system $(\mathbb{E}_2,\mathbb{M}_2)$ satisfies the condition (xi):  

\vskip+1mm
(b) the factorization systems  $(\mathbb{E}_1,\mathbb{M}_1)$ and  $(\mathbb{E}_2,\mathbb{M}_2)$ are equal;
\vskip+1mm
(c) $\mathbb{E}_1=\mathbb{E}_2$;
\vskip+1mm
(d) for any homomorphism $f:M\rightarrow M'$ with $\mathbf{R}(M'/Im\; f)=M'/Im \;f$, one has $M'=Im\;f+\mathbf{R}(M')$;
\vskip+1mm
(e) $\mathbb{M}_1=\mathbb{M}_2$;
\vskip+1mm
(f) for any module $M$ and any submodule $N$ of $M$ with  $\mathbf{R}(M)\subseteq N$, one has  $\mathbf{R}(M/N)=0$;
\vskip+1mm
(g) the factorization system $(\mathbb{E}_2,\mathbb{M}_2)$ satisfies the conditions (i), (ii) and (vi).

\end{prop}

\begin{proof}
The equivalence of the conditions (a), (b) and the one that $\mathbf{R}$ is cohereditary follows from Theorem 4.6. The equivalences  (b)$\Leftrightarrow$(c)$\Leftrightarrow$(e) follow from the property (b) of factorization systems (see the Preliminaries). The equivalences (c)$\Leftrightarrow$(d)  and (e)$\Leftrightarrow$(f) follow from the characterizations of homomorphisms from the classes  $\mathbb{E}_1$,  $\mathbb{E}_2$, $\mathbb{M}_1$, and $\mathbb{M}_2$.

If $\mathbf{R}$ is cohereditary, then, by what already was proved,  $(\mathbb{E}_1,\mathbb{M}_1)=(\mathbb{E}_2,\mathbb{M}_2)$. According to Theorem 3.1 and Proposition 3.5, the factorization system $(\mathbb{E}_1,\mathbb{M}_1)$ satisfies the conditions (i), (ii) and (vi), and hence $(\mathbb{E}_2,\mathbb{M}_2)$ also satisfies these conditions. For the converse, assume that the latter factorization system satisfies (i), (ii) and (vi). Then it is equal to $GF_1(\mathbf{R'})$, for some idempotent coherent radical $\mathbf{R'}$, as it follows from Theorem 3.1 and Proposition 3.5. Since $\mathbf{R'}$ is cohereditary, $GF_1(\mathbf{R'})=GF_2(\mathbf{R'})$, and hence $(\mathbb{E}_2,\mathbb{M}_2)=GF_2(\mathbf{R'})$. But $(\mathbb{E}_2,\mathbb{M}_2)=GF_2(\mathbf{R})$. This implies that $\mathbf{R'}=\mathbf{R}$. Therefore, $\mathbf{R}$ is cohereditary.
\end{proof}

Theorem 4.2 and Proposition 4.7 imply the following statement.

\begin{theo}
There exists a bijection between the class of cohereditary idempotent radicals and the factorization systems satisfying the conditions (i), (vii) and (xi).
\end{theo}

\vskip+2mm
\textit{Author's address:
Dali Zangurashvili, A. Razmadze Mathematical Institute of Iv. Javakhishvili Tbilisi State University},
\textit{Alexidze Str., Lane II, Tbilisi 0193, Georgia, e-mail: dali.zangurashvili@tsu.ge}

\end{document}